\documentclass{article}
\usepackage{amsmath,amssymb,amsfonts}
\usepackage{graphicx}
\usepackage{color}

\usepackage{pb-diagram}
\newtheorem{theorem}{Theorem}

\newtheorem{proposition}{Proposition}

\newtheorem{remark}{Remark}

\newenvironment{proof}[1][Proof]{\noindent\textit{#1.} }{\hfill$\Box$\medskip}

\title
      {Discriminantly separable polynomials and quad-equations}
\author{Vladimir Dragovi\'{c} and Katarina Kuki\'{c}}{}
\date{}

\begin{document}
\maketitle

\centerline{\scshape Vladimir Dragovi\'{c} }
\medskip
{\footnotesize
  \centerline{ The Department of Mathematical Sciences, University of Texas at Dallas}
   \centerline{800 West Campbell Road, Richardson TX 75080, USA}
      \centerline{and}
   \centerline{Mathematical Institute SANU}
   \centerline{Kneza Mihaila 36, 11000 Belgrade, Serbia}
}

\medskip

\centerline{\scshape Katarina Kuki\'{c}}
\medskip
{\footnotesize
  \centerline{ Faculty for Traffic and Transport Engineering, University of Belgrade}
   \centerline{Vojvode Stepe 305, 11000 Belgrade, Serbia}
}
\bigskip


\begin{abstract}
We classify the discriminantly separable polynomials of degree two
in each of three variables, defined by a  property  that all the
discriminants as polynomials of two variables are factorized as
products of two polynomials of one variable each.  Our
classification is based on the study of structures of zeros of a
polynomial component $P$ of a discriminant. This classification is
related to the classification of pencils of conics in a delicate
way. We establish a relationship between our classification and the
classification of integrable quad-equations which has been suggested
recently by Adler, Bobenko, and Suris.
\end{abstract}

\section{Introduction: An overview on discriminantly separable polynomials}\label{sec:intro}

In a  recent paper \cite {Drag3} of one of the authors of the
present paper, a new approach to the Kowalevski integration
procedure from \cite{Kow} has been suggested. It has been based on a
new notion introduced therein of {\it the discriminantly separable
polynomials}.  In next few lines we briefly give outline of those notions.

Suppose that two conics $C_1$ and $C_2$ in general position are given by
their tangential equations
\begin{equation}\label{Dragovic:equation1}
\begin{split} C_1:\,&
a_0w_1^2+a_2w_2^2+a_4w_3^2+2a_3w_2w_3+2a_5w_1w_3+2a_1w_1w_2=0;\\
 C_2:\,& w_2^2-4w_1w_3=0.
 \end{split}
 \end{equation}
We observe the pencil of conics $C(s):=C_1+sC_2$ in which conics share four common tangents.
Then the coordinate equation of the conics of the pencil is:
$$
F(s,z_1,z_2,z_3):=\det M(s, z_1,z_2,z_3)=0,
$$
with the matrix $M$ given in the next form:
\begin{equation}\label{eq:matrixM}
M(s,z_1,z_2,z_3)=\left[\begin{array}{cccc} 0 & z_1 & z_2 & z_3\\
z_1 & a_0 & a_1 & a_5 - 2s\\
z_2 & a_1 & a_2 + s & a_3\\
z_3& a_5 - 2s & a_3 & a_4
\end{array}\right].
\end{equation}
The point equation of the pencil $C(s)$ is in the form of the
quadratic polynomial in $s$
$$F:=H+Ks+Ls^2=0
$$
where $H, K,$ and $L$ are quadratic expressions in $(z_1, z_2,
z_3)$. Now we explain shortly how to introduce a new system of
coordinates in the plane, the Darboux coordinates (see \cite{Dar1}).
Given the plane with the standard coordinates $(z_1: z_2:z_3)$, we
start from the conic $C_2$ and rationally parametrize it by
$(\ell^2,\ell,1)$. The tangent line to the conic $C_2$ through a
point of the conic with the parameter $\ell_0$ is given by the
equation
$$
t_{C_2}(\ell_0): z_3\ell_0^2-2z_2\ell_0+z_1=0.
$$
Denote by $P = (\hat{z}_1, \hat{z}_2, \hat{z}_3)$ a given point in
the plane outside the conic. Each of the two solutions $x_1$ and
$x_2$ of the equation quadratic in $\ell$
$$
\hat z_3\ell^2 - 2\hat z_2\ell + \hat z_1=0
$$
corresponds to a tangent to the conic $C_2$ from the point $P$. The
pair $(x_1, x_2)$ is called the Darboux coordinates of the point
$P$. One finds immediately the converse formulae
\begin{equation}\label{eq:Darboux}
\hat z_1=x_1
x_2,\quad \hat z_2= \frac{x_1+x_2}{2},\quad  \hat z_3 =1.
\end{equation}

Changing the variables in the polynomial $F$ from the projective coordinates $(z_1:z_2:z_3)$ to the Darboux coordinates, we  rewrite its equation
$F$ in the form
\begin{equation}\label{eq:Fconics}
F(s, x_1, x_2)=L(x_1,x_2)s^2 + K(x_1,x_2)s + H(x_1,x_2)=0.
\end{equation}
In the last formula $L, K$ and $H$ are biquadratic polynomials of
$x_1,x_2$ and the explicit formulae in terms of the coefficients
 of the conic $C_1$ are:
\begin{equation}\nonumber
 \begin{split}
 L(x_1,x_2)&=(x_1-x_2)^2\\
 K(x_1,x_2)&=-a_4x_1^2x_2^2+2a_3x_1x_2(x_1+x_2)-a_5(x_1^2+x_2^2)-4a_2x_1x_2\\
 &+2a_1(x_1+x_2)-a_0\\
 H(x_1,x_2)&=(a_3^2-a_2a_4)x_1^2x_2^2+((a_1a_4-a_3a_5)x_1x_2+a_0a_3-a_5a_1)(x_1+x_2)+\\
 &(a_5^2-a_0a_4)\frac{x_1^2+x_2^2}{2}+(\frac{a_5^2-a_0a_4}{2}+2a_5a_2-2a_1a_3)x_1x_2+a_1^2-a_0a_2.
 \end{split}
\end{equation}
Here we emphasize one exceptional property that relates the equation
of pencil of conics
and the discriminants of the polynomial $F$.

The key algebraic property of the pencil equation written in the form \ref{eq:Fconics}
as a quadratic equation in
each of three variables $s, x_1, x_2$ is: {\it all three of its
discriminants are expressed as products of two polynomials in one
variable each}:
$$
\mathcal D_s(F)(x_1,x_2)=P(x_1)P(x_2),\,
\mathcal D_{x_i}(F)(s,x_j)=J(s)P(x_j), \, i,j=c.p.\{1,2\},
$$
where $J$ and $P$ are polynomials of degree $3$ and $4$
respectively. The explicit formulae for the polynomials $P$ and $J$
are
\begin{equation}\nonumber
\begin{split}
P(x)&=a_4x^4-4a_3x^3+(4a_2+2a_5)x^2-4a_1x+a_0\\
J(s)&=4s^3+4(a_2-a_5)+(a_5^2-a_0a_4+4(a_1a_3-a_2a_5))s\\
&+a_0a_3^2+a_2a_5^2++a_4a_1^2-a_0a_2a_4-2a_1a_3a_5.
\end{split}
\end{equation}
The elliptic curves
$$\Gamma_1: y^2=P(x), \quad \Gamma_2: y^2=J(s)$$
are isomorphic.

\begin{remark}\label{rem:sC_2}
The conic $C_2$ corresponds to a simple zero of the polynomial $J$.
According to our parametrization of the pencil, the value of the
parameter $s$ which corresponds to $C_2$ is $s=\infty$.
\end{remark}


In order to achieve a better understanding of the correlation
between pencils of conics and the discriminantly separable
polynomials, we briefly present a short explanation of the said
characteristics  of the pencil equation which is based on geometric
properties. A general point belongs to two conics of a tangential
pencil. If a point belongs to only one conic, then it belongs to one
of the four common tangents of the pencil and at such a point this
unique conic touches one of the four common tangents, see Fig. 1.
 \begin{figure}[h]\label{fig:2conics}
\begin{center}
\includegraphics[width=5cm,height=4cm]{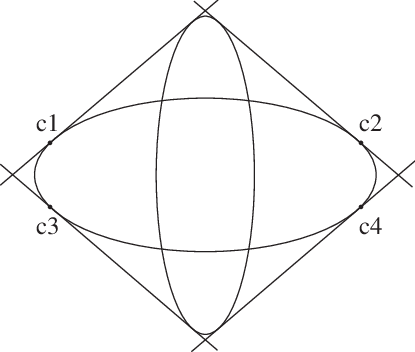}
\caption{The conics of the pencil intersect in four distinct points}
\end{center}
\end{figure}
Thus, the equation of the four common tangents may be seen as an
equation that represents the condition of annulation of the
discriminant
$$\mathcal D_s(F)(x_1,x_2)=0.$$
Previous equation is equivalent to the system
\begin{equation}\nonumber
\begin{split}
x_1&=c_1 \quad x_1=c_2 \quad x_1=c_3 \quad x_1=c_4\\
x_2&=c_1 \quad x_2=c_2 \quad x_2=c_3 \quad x_2=c_4
\end{split}
\end{equation}
where $c_i$ are parameters which correspond to the points of contact of the four common tangents with the conic $C_2$.
Finally we get
$$\mathcal D_s(F)(x_1,x_2)=P(x_1)P(x_2),$$
where the polynomial P is of the fourth degree and with four distinct roots, of the form
$$P(x)=a(x-c_1)(x-c_2)(x-c_3)(x-c_4).$$

The family of discriminantly separable polynomials in three
variables of degree two in each of them, constructed from pencils of
conics served as a motivation to introduce more general classes of
{\bf discriminantly separable polynomials}. Let us recall here the
definitions from \cite{Drag3}: a polynomial $F(x_1,\dots,x_n)$ is
{\it discriminantly separable} if there exist polynomials $f_i(x_i)$
such that for every $i=1,\dots , n$
$$
\mathcal D_{x_i}F(x_1,\dots, \hat x_i, \dots, x_n)=\prod_{j\ne
i}f_j(x_j).
$$
It is {\it symmetrically discriminantly separable} if $f_2=f_3=\dots
= f_n,$ while it is {\it strongly discriminantly separable} if
$f_1=f_2=f_3=\dots = f_n.$ It is {\it weakly discriminantly
separable} if there exist polynomials $f^i_j(x_j)$ such that for
every $i=1,\dots ,n$: $\mathcal D_{x_i}F(x_1,\dots, \hat x_i, \dots,
x_n)=\prod_{j\ne i}f^i_j(x_j). $

The so-called {\it fundamental Kowalevski equation}
\ref{eq:fundKowrel} (see \cite{Kow}, \cite{Gol})
appeared to be an example of a member of the family, as it has been
observed in \cite{Drag3}:
\begin{equation}\label{eq:fundKowrel}
Q(s, x_1,x_2):=(x_1-x_2)^2(s-\frac{l_1}{2})^2-R(x_1,x_2)(s-\frac{l_1}{2})-\frac{1}{4}R_1(x_1,x_2)=0,
\end{equation}
where $R(x_1,x_2)$ and $R_1(x_1,x_2)$ are biquadratic polynomials in $x_1$ and $x_2$ given by
\begin{equation}\nonumber
\begin{split}
R(x_1,x_2)=&-x_1^2x_2^2+6l_1x_1x_2+2lc(x_1+x_2)+c^2-k^2\\
R_{1}(x_1,x_2)=&-6l_1x_1^2x_2^2-(c^2-k^2)(x_1+x_2)^2-4lcx_1x_2(x_1+x_2)\\
&+6l_1(c^2-k^2)-4c^2l^2.
\end{split}
\end{equation}
The discriminant separability condition
$$
\mathcal D_{s}(Q)(x_1,x_2)=P(x_1)P(x_2),\,
\mathcal D_{x_i}(Q)(s,x_2)=J(s)P(x_j)
$$
is satisfied with polynomials
\begin{equation}\nonumber
\begin{split}
J(s)&=4s^3+(c^2-k^2-3l_1^2)s-l^2c^2+l_1^3-l_1k^2+l_1c^2\\
P(x_i)&=-x_i^4+6l_1 x_i^2 +4lc x_i+c^2-k^2,\, i=1,2.
\end{split}
\end{equation}
Moreover, as it has been explained in \cite{Drag3}, all the main
steps of the Kowalevski integration procedure from \cite{Kow} (see
also \cite{Gol}) now follow as easy and transparent
logical consequences of the theory of discriminantly separable
polynomials.

There are two natural and important questions in this context:

1) {\it Are there any other discriminantly separable polynomials of
degree two in each of three variables, beside those constructed from
pencils of conics? In addition, the question is to perform a
classification of such polynomials}.

2) {\it Are there other integrable dynamical systems related to
discriminantly separable polynomials?}


{\bf The main issue of this paper is to address these two key
questions.}


In order to make precise the first question, one needs to specify a
gauge group or the classes of equivalence up to which  a
classification would be performed. This leads to the group of
M\"obius transformations, as introduced in Corollary 3 of
\cite{Drag3}:
\begin{equation}\label{eq:mobius}
 x_1\mapsto \frac {a_1x_1+b_1}{c_1x_1+d_1},\,
x_2\mapsto \frac {a_2x_2+b_2}{c_2x_2+d_2},\, s\mapsto \frac
{a_3s+b_3}{c_3s+d_3}.
\end{equation}
The classification of strongly discriminantly separable polynomials
$\mathcal{F}(x_1,x_2,x_3)$ of degree two in each of three variables
modulo fractional-linear transformations from formulae
\ref{eq:mobius} with $a_1=a_2=a_3,\, b_1=b_2=b_3,\,
c_1=c_2=c_3,\,d_1=d_2=d_3, $ as the gauge group, is one of the main
tasks of the present paper.

There is a remarkable correspondence between this classification and the classification of pencils of conics,
see \cite{DragRad}, for example, for more details
about pencils of conics. In the case of general position, the
conics of a pencil intersect in four distinct points, and we code
such situation with $(1,1,1,1)$, see Fig. 2, and denote it by (A). It
corresponds to the case where the polynomial $P$ has four simple zeros.
In this case, the family of strongly discriminantly separable
polynomials coincides with the family constructed above
from a general pencil of conics. This family, as it has been
indicated in \cite{Drag3}, corresponds to the two-valued
Buchstaber-Novikov group associated with a cubic curve
$$
\Gamma_2: y^2=J(s).
$$
Two-valued groups of the form $(\Gamma_2, Z_2)$ have been introduced
in \cite{Buc1}. The so-called Kowalevski change of variables
appeared to be an infinitesimal of the two-valued operation in this
group, see \cite{Drag3}. (The theory of $n$ valued-groups originates
from a pioneering paper \cite{BN}. For a modern account see
\cite{Buc}, and for higher-genus analogues, see \cite{BD}.)

However, in the degenerate cases the above mentioned correspondence
between the discriminantly separable polynomials and pencils of
conics is much more delicate. It remains valid in all cases where
the polynomial $J$ has at least one simple zero (such cases are
coded as $(1,1,2)$ and $(1,3)$). Surprisingly enough, the
correspondence is broken in the other two cases, which are
characterized by the fact that the polynomial $J$ has  multiple
zeros only, and which are coded as $(2,2)$ and $(4)$.

The polynomials which appear in the case of tangential equation of a
general pencil of conics as well as in the fundamental Kowalevski
equation \ref{eq:fundKowrel} are symmetrically discriminantly
separable. Let us remark that by a suitable M\"{o}bius
transformation it can be transformed into an equivalent strongly
discriminantly separable polynomial, see \cite{Drag3}.


Referring to the second key question, one direction related to the
continuous systems has been derived recently in \cite{DK} and
\cite{DK1}. Here we trace a connection between the discriminantly
separable polynomials of degree two in each of three variables with
the quad-equations. See \cite {ABS1} and \cite {ABS2} for more
details about the quad-equations; some basic notions from there are
collected in Section \ref{sec:ABS}. We establish a correspondence
between the discriminantly separable polynomials of degree two in
three variables and the families of symmetric biquadratic
polynomials with discriminants independent on parameter. Latter
appear in the construction of Adler, Bobenko, and Suris \cite{ABS2}.
This correspondence enables us to assign certain quad-equations to
the discriminantly separable polynomials.

Recently, a new class of geometric quad-graphs, of the so-called
line congruences, associated with pencils of quadrics has been
introduced in \cite{DragRad1}.


\section{Classification of strongly discriminantly separable polynomials of degree two in three
variables}\label{sec:class}

In this section we will classify strongly discriminantly separable
 polynomials
$\mathcal{F}(x_1,x_2,x_3)\in\mathbb{C}[x_1,x_2,x_3]$ which are of
degree two in each variable, modulo a group of the M\"{o}bius transformations
\begin{equation}\label{eq:mobius2}
x_1\mapsto \frac {ax_1+b}{cx_1+d},\, x_2\mapsto \frac
{ax_2+b}{cx_2+d},\, x_3\mapsto \frac {ax_3+b}{cx_3+d}.
\end{equation}

Let
\begin{equation}\label{eq:F(x_1,x_2,x_3)}
\mathcal{F}(x_1,x_2,x_3)=\sum_{i,j,k=0}^{2}a_{ijk}x_{1}^{i}x_{2}^{j}x_{3}^{k}
\end{equation}
be a strongly discriminantly separable polynomial with
\begin{equation}\label{eq:SDS}
\mathcal{D}_{x_i}\mathcal{F}(x_j,x_k)=P(x_j)P(x_k),\quad (i,j,k)=c.p.(1,2,3).
\end{equation}

Here, by $\mathcal{D}_{x_{i}}\mathcal{F}(x_{j},x_{k})$ we denote the
discriminant of $\mathcal{F}$ considered as a quadratic polynomial
in $x_i$.


By plugging  \ref{eq:F(x_1,x_2,x_3)} into \ref{eq:SDS} for a
given polynomial $P(x)=Ax^4+Bx^3+Cx^2+Dx+E,$ we get a system of $75$
equations of degree two with $27$ unknowns $a_{ijk}$.


\begin{theorem}\label{th:classification}
The strongly discriminantly separable polynomials
$\mathcal{F}(x_1,x_2,x_3)$ satisfying \ref{eq:SDS} modulo
fractional linear transformations are exhausted by the following
list depending on the structure  of the roots of a non-zero
polynomial $P(x)$ :
\begin{description}

    \item[(A)] $P$ has four simple zeros, with the canonical form
    $P_A(x)=(k^2x^2-1)(x^2-1)$,
     \begin{equation}\nonumber
     \begin{split}
\mathcal{F}_A=&\frac{1}{2}(-k^2x_1^2-k^2x_2^2+1+k^2x_1^2x_2^2)x_3^2+(1-k^2)x_1x_2x_3\\
&+\frac{1}{2}(x_1^2+x_2^2-k^2x_1^2x_2^2-1),
\end{split}
\end{equation}

   \item[(B)] $P$ has two simple zeros and one double zero, with the canonical form
    $P_B(x)=x^2-e^2,\,e \neq 0$,
     $$ \mathcal{F}_B=x_1x_2x_3 +\frac{e}{2}(x_1^2+x_2^2+x_3^2-e^2),$$

    \item[(C)] $P$ has two double zeros, with the canonical form $P_C(x)=x^2$,
$$\mathcal{F}_{C_1}=\lambda x_1^2x_2^2+\mu x_1x_2x_3+\nu x_3^2,\quad
    \mu^2-4\lambda\nu=1,$$
    $$\mathcal{F}_{C_2}=\lambda x_1^2x_3^2+\mu x_1x_2x_3+\nu x_2^2,\quad
    \mu^2-4\lambda\nu=1,$$
    $$\mathcal{F}_{C_3}=\lambda x_2^2x_3^2+\mu x_1x_2x_3+\nu x_1^2,\quad
    \mu^2-4\lambda\nu=1,$$
         $$\mathcal{F}_{C_4}=\lambda x_1^2x_2^2x_3^2+\mu x_1x_2x_3+\nu,\quad
    \mu^2-4\lambda\nu=1,$$

 \item[(D)] $P$ has one simple and one triple zero,with the canonical form $P_D(x)=x$,
 \begin{equation}\nonumber \mathcal{F}_D=-\frac{1}{2}(x_1x_2+x_2x_3+x_1x_3)+\frac{1}{4}(x_1^2+x_2^2+x_3^2),\end{equation}

 \item[(E)] $P$ has one quadruple zero, with the canonical form $P_E(x)=1$,
    $$\mathcal{F}_{E_1}= \lambda(x_1+x_2+x_3)^2+\mu(x_1+x_2+x_3)+\nu,\quad \mu^2-4\lambda\nu=1,$$
    $$\mathcal{F}_{E_2}=
\lambda(x_2+x_3-x_1)^2+\mu(x_2+x_3-x_1)+\nu,\quad
\mu^2-4\lambda\nu=1,$$
    $$\mathcal{F}_{E_3}=
\lambda(x_1+x_3-x_2)^2+\mu(x_1+x_3-x_2)+\nu,\quad
\mu^2-4\lambda\nu=1,$$
    $$\mathcal{F}_{E_4}=
\lambda(x_1+x_2-x_3)^2+\mu(x_1+x_2-x_3)+\nu,\quad
\mu^2-4\lambda\nu=1.$$

\end{description}
\end{theorem}

\begin{proof}
Proof is done by a straightforward calculation by solving the system
of equations \ref{eq:SDS} for the canonical representatives of the
polynomials P. For example, in the case when
$$P_A(x)=(k^2x^2-1)(x^2-1),$$ as the solutions of the system \ref{eq:SDS}, up to the sign,
we get the polynomials:
\begin{equation}\nonumber
\begin{split} \mathcal{F}_{A_1} = &(-\frac{k^2}{2}x_1^2-\frac{k^2}{2}x_2^2+\frac{1}{2}+
\frac{k^2}{2}x_1^2x_2^2)x_3^2+(1-k^2)x_1x_2x_3\\
&+
\frac{1}{2}x_1^2+\frac{1}{2}x_2^2-\frac{k^2}{2}x_1^2x_2^2-\frac{1}{2},\\
\mathcal{F}_{A_2} = &(-\frac{k^2}{2}x_1^2-\frac{k^2}{2}x_2^2+\frac{1}{2}+\frac{k^2}{2}x_1^2x_2^2)x_3^2+(-1+k^2)x_1x_2x_3\\
&+
\frac{1}{2}x_1^2+\frac{1}{2}x_2^2-\frac{k^2}{2}x_1^2x_2^2-\frac{1}{2},
\end{split}
\end{equation}
\begin{equation}\nonumber
\begin{split}
\mathcal{F}_{A_3} = & (\frac{k}{2}-\frac{k}{2}x_1^2-\frac{k}{2}x_2^2+\frac{k^3}{2}x_1^2x_2^2)x_3^2+(1-k^2)x_1x_2x_3\\
&+\frac{k}{2}x_1^2+\frac{k}{2}x_2^2-\frac{k}{2}x_1^2x_2^2-\frac{1}{2k},\\
\mathcal{F}_{A_4} = & (\frac{k}{2}-\frac{k}{2}x_1^2-\frac{k}{2}x_2^2+\frac{k^3}{2}x_1^2x_2^2)x_3^2+(-1+k^2)x_1x_2x_3\\
&+\frac{k}{2}x_1^2+\frac{k}{2}x_2^2-\frac{k}{2}x_1^2x_2^2-\frac{1}{2k}.
\end{split}
\end{equation}

The polynomials $\mathcal{F}_{A_1}$ and
$\mathcal{F}_{A_2}$ are  equivalent modulo the
fractional-linear transformations \ref{eq:mobius2} $x_i \mapsto
-x_i, \quad i=1,2,3$. The same is true for the polynomials $\mathcal{F}_{A_3}$ and $\mathcal{F}_{A_4}$.
Further, acting by the transformations $x_i \mapsto 1/(kx_i),\,
i=1,2,3$ on $\mathcal{F}_{A_3}(x_1,x_2,x_3)$ we get the gauge
equivalence between $\mathcal{F}_{A_1}$ and
$\mathcal{F}_{A_3}$. In the same way we get the polynomials in the cases (B)-(E).
\end{proof}

The correspondence between this classification and pencil of conics
in the case (A) has been briefly presented in Introduction. The
corresponding pencil of conics is presented on Fig. 2.

\vspace{-.5cm}
\begin{figure}[h]\label{fig:4points}
    \begin{center}
    \includegraphics[width=3.5cm,height=3.5cm]{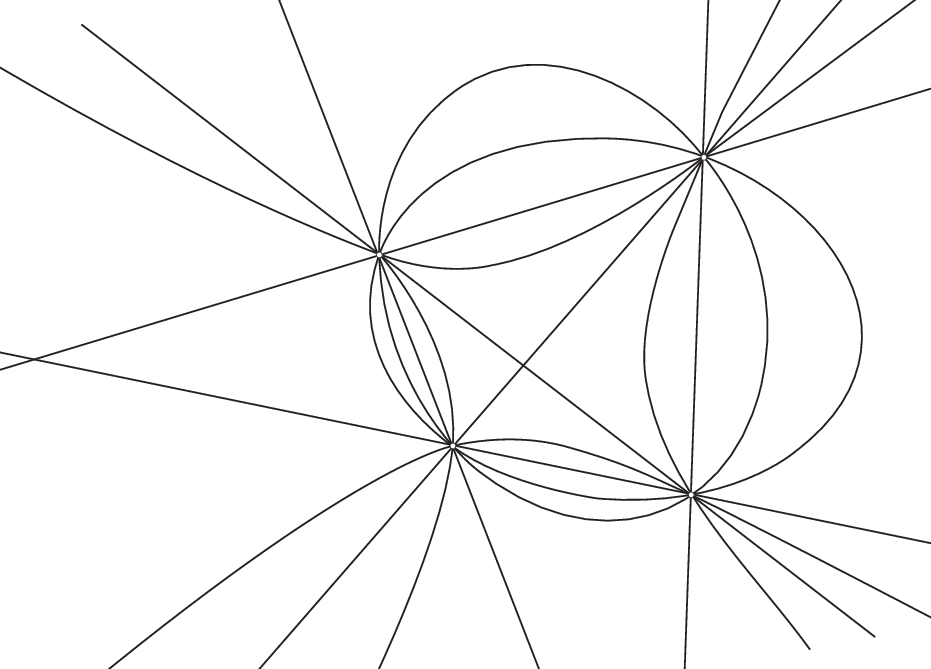}
    \caption{Pencil with four simple points}
    \end{center}
    \end{figure}
\vspace{-.3cm}

Without losing generality, we will keep using the conic $C_2$ as one
with respect to which we define the Darboux coordinates. The
obtained families of polynomials in cases (A), (B), and (D) are
unique up to M\"obius equivalence  (each being M\"{o}bius-equivalent
to a corresponding strongly discriminantly separable polynomial),
and they  represent the equations of pencils of conics of the type
(A), (B), and (D). The situation in the cases (C) and (E) is
dramatically different. Not only are we losing uniqueness, up to
M\"{o}bius equivalence, of families of the polynomials, but also we
are losing such a transparent geometric
correlation with pencils of conics. We will present now the connection with pencils of conics in the cases (B) and (D). Later on, we will discuss
the cases (C) and (E).
%

We consider the case (B) coded $(1,1,2)$ when the polynomial $P$ has
two simple zeros and one double zero. The conics of the
corresponding pencil intersect in two simple points, and they have a
common tangent in the third point of intersection. We start from
conics $C_1$ and $C_2$ given by \ref{Dragovic:equation1}. In this
case, we search for the conic $C_1$ which has the intersection with
$C_2$ containing two simple points $Q$ and $S$ and at the third
point of intersection $R$, the conics have a common tangent line.
Without loosing generality, choose points from $C_2$: $R=(1,1,1),
 Q=(0,0,1)$ and $S=(a^2,a,1)$, with $a \neq 0,1$ in an appropriate projective frame. The equation of the common tangent at $R$ is $t:z_1-2z_2+z_3=0$.

We get the following coefficients of $C_1$ in the tangential equation:
\begin{equation}\label{eq:coeffB}
\begin{split}
a_2&=\frac{a_1^2}{a_0}\\
a_3&=\frac{aa_1+2a_1-2a_0}{a}\\
a_4&=\frac{aa_0+4a_1-4a_0}{a}\\
a_5&=-\frac{aa_0^2-4aa_0a_1+2aa_1^2-2a_0a_1+2a_0^2}{aa_0},
\end{split}
\end{equation}
with an arbitrary $a_1$ and $a_0\neq 0$. Then, the coordinate
equation of the conics of the pencil $C_1+sC_2$ is
$$F(s,z_1,z_2,z_3):=\det M(s,z_1,z_2,z_3)=0.$$
The matrix $M$ is obtained by plugging the coefficients
\ref{eq:coeffB} into \ref{eq:matrixM}. Finally in the Darboux
coordinates $(x_1,x_2)$ as in \ref{eq:Darboux} we get the equation
of the pencil $C(s)$ in the following form:
$$F(x_1,x_2,s):=L(x_1,x_2)s^2-\frac{K(x_1,x_2)}{aa_0}s+\frac{(a_0-aa_1)(a_1-a_0)^2H(x_1,x_2)}{a_0^2a^2}=0$$
with
\begin{equation}\nonumber
\begin{split}
L(x_1,x_2)&=(x_1-x_2)^2\\
K(x_1,x_2)&=a_0(aa_0+4a_1-4a_0)x_1^2x_2^2-2a_0(aa_1+2a_1-2a_0)x_1x_2(x_1+x_2)\\
&+(2a_1a_0(1+2a)-aa_0^2-2aa_1^2-2a_0^2)(x_1^2+x_2^2)+4aa_1^2x_1x_2\\
&-2aa_0a_1(x_1+x_2)+aa_0^2\\
H(x_1,x_2)&=4a_0x_1^2x_2^2-2a_0(a+2)x_1x_2(x_1+x_2)+(a_0+2aa_0-aa_1)(x_1^2+x_2^2)\\
&+2(a_0+2aa_0+aa_1)x_1x_2-2aa_0(x_1+x_2).
\end{split}
\end{equation}
The discriminants of $F$ are factorized as products of polynomials
with two simple and one double zero:
\begin{equation}\nonumber
\begin{split}
P(x)&=(x-1)^2((4(a_0-a_1)-aa_0)x^2+2ax(2a_1-a_0)-aa_0)\\
J(s)&=(sa_0+(a_0-a_1)^2)(saa_0+(aa_1-a_0)(a_1-a_0))^2.
\end{split}
\end{equation}
\begin{figure}[h]\label{fig:31points}
     \begin{center}
\includegraphics[width=3.5cm,height=3.5cm]{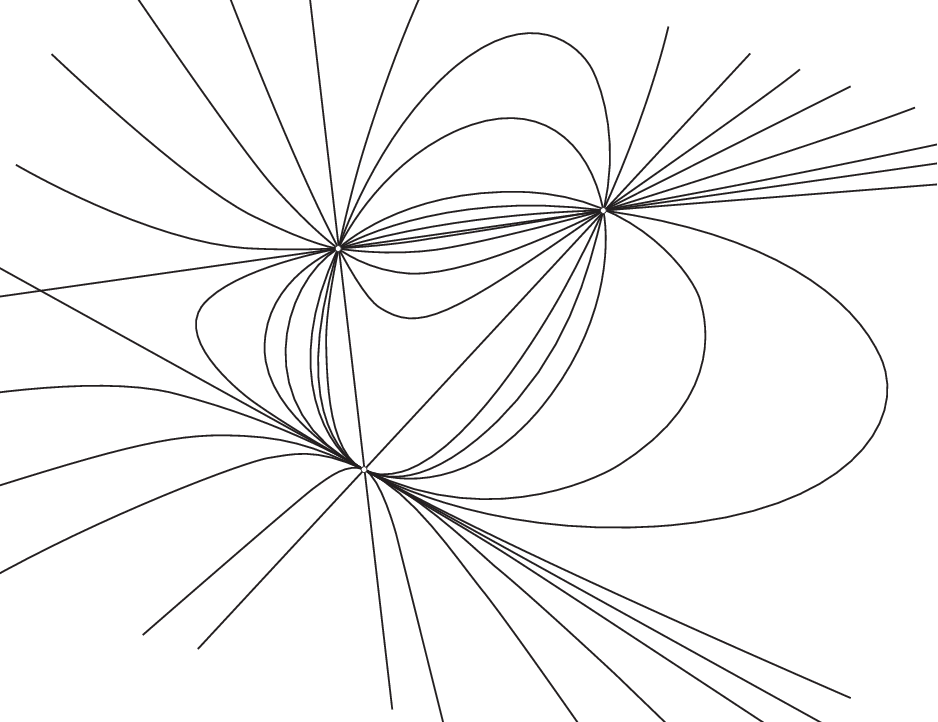}
\caption{Pencil with one double and two simple points}
\end{center}
\end{figure}
\begin{remark} The roots of the polynomial $J$ are allocated to degenerate conics from the
pencil $C(s)$. For the double root of the polynomial $J$, we get the
conic
$$C\left(\frac{(a_0-a_1)(aa_1-a_0)}{aa_0}\right):=(x_1-1)^2(x_2-1)^2=0$$
which is the  double tangent line at $R$. Similarly, for the simple
root of the polynomial $J$ we get the double line
$$C\left(-\frac{(a_1-a_0)^2}{a_0}\right):=(x_1+x_2+(a-2)x_1x_2-a)^2=0.$$
\end{remark}

As in the Remark \ref{rem:sC_2}, the conic $C_2$ corresponds  to $s=\infty$.

Now we consider case (D). In this case the polynomial $P$ has one
simple zero and one triple zero, coded $(1,3)$. The conics of the
corresponding pencil intersect at one simple point $S$, and they
have another  common point $R$, with the common tangent line and the
same curvature at that point. As in the case (B) we choose points
from $C_2$: $R=(1,1,1)$ and $S=(a^2,a,1)$ in an appropriate
projective frame. The common tangent line at $R$ is
$t:z_1-2z_2+z_3=0$ and the curvature of $C_1$ at $R$ is equal to the
curvature of $C_2$ at that point, which is $\displaystyle
\frac{5\sqrt{5}}{2}$. Solving that system of equations we get
coefficients of the tangential equation of $C_1$:
\begin{equation}\label{eq:coeffD}
\begin{split}
a_0&=-3a_4+4a_3\\
a_1&=-2a_4+3a_3\\
a_2&=\frac{12a_3a_4+aa_3^2-4a_4^2-9a_3^2}{3a_4-4a_3+aa_4}\\
a_5&=\frac{6a_3a_4(1+a)-a_4^2-6a_3^2-3aa_4^2-2aa_3^2}{3a_4-4a_3+aa_4}.
\end{split}
\end{equation}
After plugging these coefficients into the matrix \ref{eq:matrixM}
and introducing the Darboux coordinates \ref{eq:Darboux} we get
the equation of the pencil
$$F(x_1,x_2,s):=L(x_1,x_2)s^2-\frac{K(x_1,x_2)}{3a_4-4a_3+aa_4}s+\frac{(a_4-a_3)^3H(x_1,x_2)}{(3a_4-4a_3+aa_4)^2}=0$$
with
\begin{equation}\nonumber
\begin{split}
L(x_1,x_2)&=(x_1-x_2)^2\\
K(x_1,x_2)&=(aa_4^2-4a_3a_4+3a_4^2)x_1^2x_2^2-2a_3(3a_4-4a_3+aa_4)x_1x_2(x_1+x_2)\\
&+(6a_3a_4(1+a)-2aa_3^2-3aa_4^2-a_4^2-6a_3^2)(x_1^2+x_2^2)\\
&+(48a_3a_4+4aa_3^2-36a_3^2-16a_4^2)x_1x_2\\
&+2(2a_4-3a_3)(3a_4-4a_3+aa_4)(x_1+x_2)\\
&-(3a_4-4a_3)(3a_4-4a_3+aa_4)\\
H(x_1,x_2)&=4(aa_4+3a_4-4a_3)x_1^2x_2^2-2(a+3)(3a_4-4a_3+aa_4)x_1x_2(x_1+x_2)\\
&+(3a_4a^2+6aa_4+7a_4-a^2a_3-6aa_3-9a_3)(x_1^2+x_2^2)\\
&+(22a_4-36aa_3+6a_4a^2-30a_3+2a^2a_3+36aa_4)x_1x_2\\
&-2(1+3a)(3a_4-4a_3+aa_4)(x_1+x_2)+4a(3a_4-4a_3+aa_4).
\end{split}
\end{equation}
The discriminants of $F$ are factorized as products of polynomials
with one simple and one triple zero:
\begin{equation}\nonumber
\begin{split}
P(x)&=(x-1)^3(a_4x_2+3a_4-4a_3)\\
J(s)&=((3a_4-4a_3+aa_4)s+(a-1)(a_3-a_4)^2)^3.
\end{split}
\end{equation}
\begin{remark}
Plugging the triple root of the polynomial $J$ into the equation of
the pencil, we get a degenerative conic, the double common tangent
line
$$C\left(\frac{(1-a)(a_3-a_4)^2}{3a_4-4a_3+aa_4}\right):=(x_1-1)^2(x_2-1)^2=0.$$
\end{remark}

\begin{figure}[h]\label{fig:13points}
    \begin{center}
    \includegraphics[width=3.5cm,height=3.5cm]{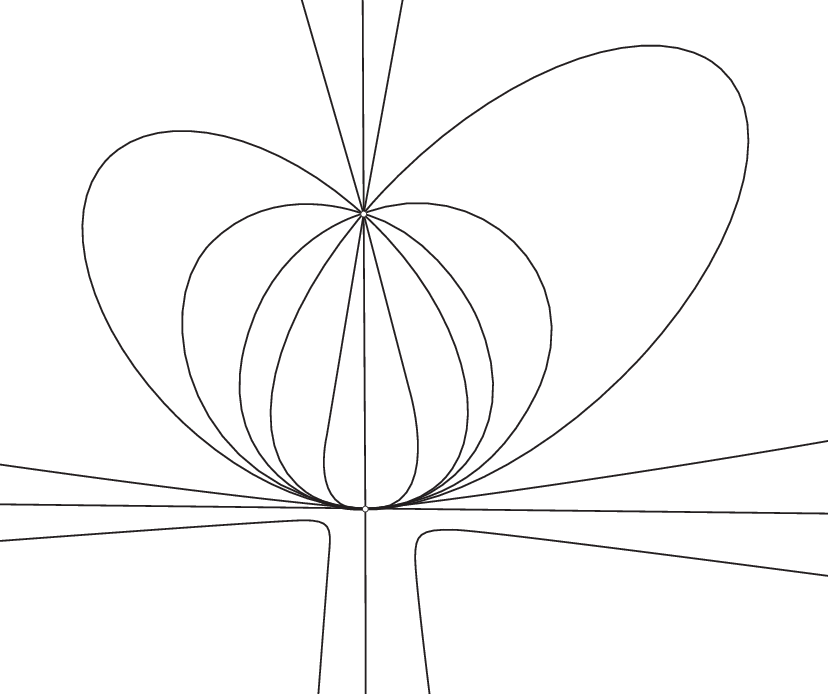}
    \caption{Pencil with one simple and one triple point}
    \end{center}
    \end{figure}

From the Remark \ref{rem:sC_2}, we know that the polynomial $J$ has
a simple zero that corresponds to the conic $C_2$. Regarding the
cases (C) and (E) this shows that starting from the corresponding
pencils of conics, it is not possible to get $J$ which has  multiple
zeros only. Indeed, the next simple calculation validates the claim.

For  the case $(2,2)$ we search for a conic $C_1$ that intersects
with $C_2$ in two points $R$ and $S$, and which has a common tangent
with $C_2$  in each of those points. In an appropriate projective
frame, we choose $R(1,1,1)$ and $S(a^2,a,1)$ with $a\neq 1$, and
then, the common tangents  are $t_1:z_1-2z_2+z_3=0$ and
$t_2:z_1-2az_2+a^2z_3=0.$ The coefficients of the tangential
equation of $C_1$ are
\begin{equation}\label{eq:coeffC}
\begin{split}
a_0&=a^2a_4\\
a_1&=\frac{aa_4(a+1)}{2}\\
a_2&=\frac{a^2a_4+4aa_4-2a_5+a_4}{4}\\
a_3&=\frac{(a+1)a_4}{2}.
\end{split}
\end{equation}
After plugging \ref{eq:coeffC} into the matrix \ref{eq:matrixM}
and introducing the Darboux coordinates \ref{eq:Darboux}, we get
the equation of a pencil:
$$F(x_1,x_2,s):=L(x_1,x_2)s^2+K(x_1,x_2)s+\frac{(aa_4-a_5)H(x_1,x_2)}{4}=0$$
with
\begin{equation}\nonumber
\begin{split}
L(x_1,x_2)&=(x_1-x_2)^2\\
K(x_1,x_2)&=-a_4x_1^2x_2^2+a_4(a+1)x_1x_2(x_1+x_2)-a_5(x_1^2+x_2^2)\\
&+(2a_5-a_4-4aa_4-a^2a_4)x_1x_2+aa_4(a+1)(x_1+x_2)-a^2a_4\\
H(x_1,x_2)&=-2a_4x_1^2x_2^2+2a_4(a+1)x_1x_2(x_1+x_2)-(aa_4+a_5)(x_1^2+x_2^2)\\
&+(2a_5-2a_4-2a^2a_4-6aa_4)x_1x_2+2aa_4(a+1)(x_1+x_2)-2a^2a_4.
\end{split}
\end{equation}
Notice here that $F$ can be factorized as
\begin{equation}\nonumber
\begin{split}
F(x_1,x_2,s)&=\left(\frac{2s-a_5+aa_4}{4}\right)(2(x_1-x_2)^2s-2a_4x_1^2x_2^2\\
&+2a_4(a+1)x_1x_2(x_1+x_2)+(2a_5-2a_4-2a^2a_4-6aa_4)x_1x_2\\
&-(aa_4+a_5)(x_1^2+x_2^2)+2aa_4(a+1)(x_1+x_2)-2a^2a_4).
\end{split}
\end{equation}
The discriminants of $F$ are factorized as products of the following
polynomials:
\begin{equation}\nonumber
\begin{split}
P(x)&=(x-1)^2(x-a)^2\\
J(s)&=(2s+aa_4-a_5)^2(4s-2a_5+a_4(1+a^2)).
\end{split}
\end{equation}
\begin{remark}
Plugging a simple root of the polynomial $J$ into the equation of
the pencil we get the double line
$$C\left(\frac{2a_5-a_4(1+a^2)}{4}\right):=((x_1+x_2)(1+a)-2x_1x_2-2a)^2=0$$
But, plugging the double root $\displaystyle s_0=\frac{a_5-aa_4}{2}$
we get $F(x_1,x_2,s_0)\equiv0$.
\end{remark}

\begin{figure}[h]\label{fig:22points}
    \begin{center}
    \includegraphics[width=3.5cm,height=3.5cm]{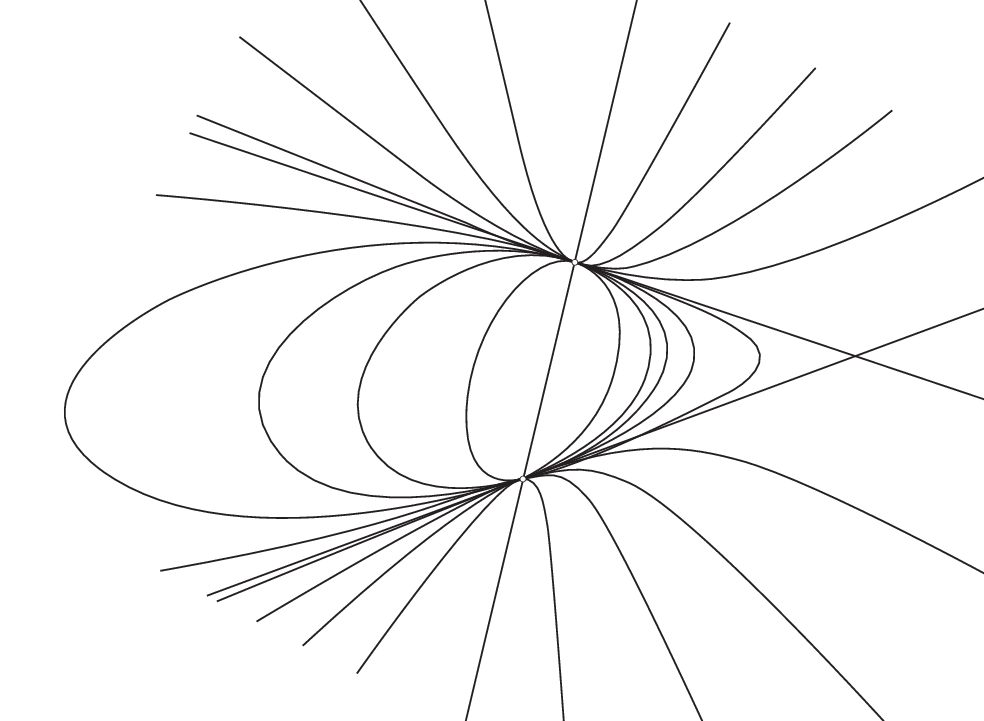}
    \caption{Pencil with two double points}
    \end{center}
    \end{figure}

This unexpected lack of  corresponding pencils of conics in the
cases (C) and (E) can be understood better in the light of the
following statement:

\begin{theorem}[Corollary, \cite{SK}, VIII, Ch. 1]
The symmetrical $(2-2)$ algebraic correspondence cut on $C$ by
tangents to $C_2$ breaks up into a homography and its inverse if and
only if $C_2$ has double contact with $C$. If $C_2$ coincides with
$C$, the correspondence is the identical correspondence taken twice.
\end{theorem}

The two points of contact correspond to the fixed point of the
homography (a M\"obius transformation in our terminology).

If we fix a value of the parameter $s$ in the equation of the pencil
$F(x_1,x_2,s)=0$, then
$$F(x_1, x_2, s)=\hat F_s(x_1,x_2)=0$$
specifies a $(2-2)$ correspondence $\hat F_s$. Here $x_1$ is the
parameter on the conic $C$ that corresponds to a point on $C$, and
$x_2$ corresponds to another point of $C$. We also suppose that the
two fixed points of the M\"obius transformation (homography) are the
points with the parameters equal to 0 and $\infty$. The M\"{o}bius
transformation, denoted  $w$, is then of the form $w(x)=ax$. In a
case of the pencil of conics $(2, 2)$ with the intersection at two
double points, according to the previous Theorem, the polynomial
$\hat F_s(x_1, x_2)$ must be of the form
\begin{equation}\label{eq:FCTheorem}
\hat F_s(x_1,x_2)=(ax_1+bx_2)(bx_1+ax_2).
\end{equation}
However, one can easily see that for a fixed value of the parameter
$s$, the polynomials $\mathcal{F}_{C_1}-\mathcal{F}_{C_4}$, do not
have the form \ref{eq:FCTheorem}. This shows that those
polynomials do not correspond to pencils of conics with two double
base points. The same consideration  when the fixed points of $w$
coincide, explains the case (E).

\begin{figure}[h]\label{fig:14points}
    \begin{center}
    \includegraphics[width=3.5cm,height=3.5cm]{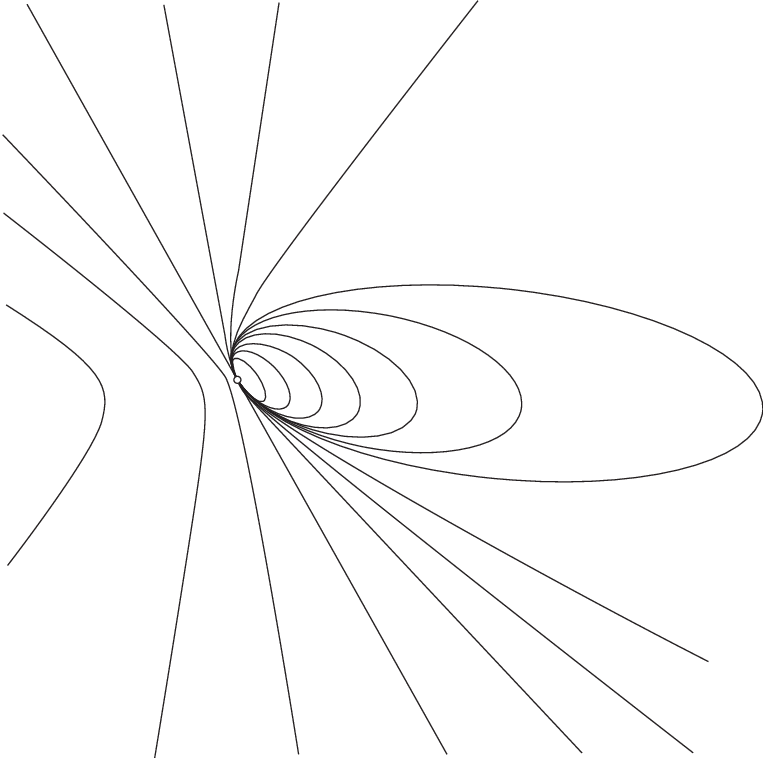}
    \caption{Pencil with one quadruple point}
     \end{center}
    \end{figure}


\section{From discriminantly separable polynomials to integrable quad-equat\-ions}\label{sec:ABS}

Now, we will emphasize a connection between the discriminantly
separable polynomials and integrable quad-equations. Latter has been
developed by Adler, Bobenko, Suris \cite{ABS1}, \cite{ABS2}, see
also \cite{BS1}, \cite{BS2}.


Recall that the quad-equations
are the equations on quadrilaterals of the form
\begin{equation}\label{eq:quadeq}
Q(x_1,x_2,x_3,x_4)=0
\end{equation}
where $Q$ is a polynomial of degree one in each argument, i.e. $Q$
is a multiaffine polynomial. The field variables $x_i$ are assigned to four
vertices of a quadrilateral as in a Figure 6. Besides depending on the variables $x_1,...,x_4\in \mathbb{C}$, the polynomial Q
still depends on two parameters $\alpha,\beta\in \mathbb{C}$ that are assigned to the edges of a quadrilateral. It is assumed that opposite edges carrie the same parameter.

\begin{figure}[h]\label{fig:quadrilateral}
    \begin{center}
    \includegraphics[width=3cm,height=3cm]{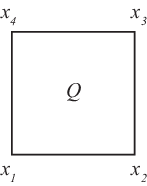}
    \caption{An elementary quadrilateral.}
     \end{center}
    \end{figure}


The equation \ref{eq:quadeq} can be solved for each variable, and the solution
is a rational function of the other three variables.  A solution
$(x_1,x_2,x_3,x_4)$ of equation \ref{eq:quadeq} is
\textit{singular} with respect to $x_i$ if it also satisfies the
equation $Q_{x_i}(x_1,x_2,x_3,x_4)=0$.


Following \cite{ABS2} we consider the idea of integrability as {\it
a consistency}, see Figure 8. We assign six quad-equations to the
faces of the coordinate cube. The system is said to be
\textit{3D-consistent} if the three values for $x_{123}$ obtained
from the equations on the right, back and top faces coincide for
arbitrary initial data $x,\,x_1,\,x_2,\,x_3$.

\begin{figure}[h]\label{fig:cube}
    \begin{center}
    \includegraphics[width=3.5cm,height=3.5cm]{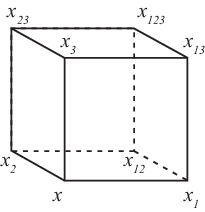}
    \caption{A 3D consistency. }
    \end{center}
    \end{figure}

Then, by applying the discriminant-like operators introduced in
\cite{ABS2}
\begin{equation}\label{eq:discrop}
\delta_{x,y}(Q)=Q_xQ_y-QQ_{xy},\quad \delta_x(h)=h_x^2-2hh_{xx},
\end{equation}
one can make descent from the faces to the edges and then to the vertices of the cube: from a  multiaffine polynomial
$Q(x_1,x_2,x_3,x_4)$ to  a
biquadratic polynomial $h(x_i,x_j):=\delta_{x_k,x_l}(Q(x_i, x_j,
x_k, x_l))$ and further to a polynomial
$P(x_i)=\delta_{x_j}(h(x_i,x_j))$ of degree up to four.
 By  using  of relative invariants of polynomials under fractional linear transformations, the formulae that
express $Q$ through  the biquadratic polynomials of three edges,
were derived in \cite{ABS2}:
\begin{equation}\label{eq:QABS2}\frac{2Q_{x_1}}{Q}=\frac{h^{12}_{x_1}h^{34}-h^{14}_{x_1}h^{23}+h^{23}h_{x_3}^{34}-h_{x_3}^{23}h^{34}}{h^{12}h^{34}-h^{14}h^{23}}.
\end{equation}


 A biquadratic polynomial $h(x,y)$ is said to be \textit{nondegenerate} if no polynomial in its equivalence class
 with respect to the fractional linear transformations, is divisible by a factor of the form $x-c$ or $y-c$, with $c=const$.
 A multiaffine function $Q(x_1,x_2,x_3,x_4)$ is said to be of \textit {type $Q$} if all four of its accompanying
 biquadratic polynomials $h^{jk}$ are nondegenerate.  Otherwise, it is of \textit{type $H$}. Previous
notions were introduced in \cite{ABS2}.


Now, we start with an arbitrary strongly discriminantly separable
polynomial $
\mathcal{F}(x_1,x_2,\alpha)$ of degree two in each of three variables. In order to relate that polynomial with corresponding
quad-equations we have to provide corresponding biquadratic polynomial
$h=h(x_1, x_2)$ and a multiaffine polynomial $Q=Q(x_1, x_2, x_3,\break
x_4)$.

The requirement that the discriminants of $h(x_1,x_2)$ do not depend on $\alpha$,
 see \cite{ABS1}, \cite{ABS2}, will be satisfied if as a biquadratic polynomials
  $h(x_1,x_2)$  we take
  $$\hat h(x_1, x_2):=\displaystyle
  \frac{\mathcal{F}(x_1,x_2,\alpha)}{\sqrt{P(\alpha)}}.$$

\begin{proposition}\label{proph}
The biquadratic polynomials \begin{equation}\label{eq:corrhF}
\hat{h}_{I}(x_1,x_2)=\frac{\mathcal{F}_{I}(x_1,x_2,\alpha)}{\sqrt{P_{I}(\alpha)}}
\end{equation}
 satisfy
$$\delta_{x_1}(\hat h)=P_{I}(x_2),\quad \delta_{x_2}(\hat h)=P_{I}(x_1)$$
for $I=A,B,C,D,E$ and polynomials $P_I, \mathcal{F}_{I}$ from Theorem 1.
\end{proposition}

Now from the formulae \ref{eq:QABS2} and with the polynomials $h^{ij}$ replaced by
$\hat{h}^{ij}$ we can easily get quad-equations which correspond to representatives of discriminantly separable polynomials from
Theorem \ref{th:classification}. Those equations represent re-parametrization of quad-equations from the list obtained in \cite{ABS2}.

The important feature of the quad-equations obtained from the biquadratic polynomials $\hat{h}(x_1,x_2)$ is that parameter $\alpha$ plays the same role as $x_1$
and $x_2$.

\section*{Acknowledgments}
The authors would like to thank to Prof. Yu. B. Suris for his kind
advices and very instructive remarks. The research was partially
supported by the Serbian Ministry of Science and Technological
Development, Project 174020 {\it Geometry and Topology of Manifolds,
Classical Mechanics and Integrable Dynamical Systems}. The authors
are grateful to the referees for comments and suggestions.

\end{document}